\documentclass[12pt]{article}
\usepackage{hyperref,amsfonts,amssymb}
\usepackage{amsfonts,amssymb}
\usepackage{mathtools}
\usepackage{rotating}

\def\C{\mathbb{C}}
\def\A{\mathbb{A}}

\def\bq{ \begin{equation} }
\def\eq{ \end{equation} }
\def\ben{ \begin{eqnarray} }
\def\en{ \end{eqnarray} }

\def\frac#1#2{{#1\over #2}}

\def\on#1#2{\mathop{\vbox{\ialign{##\crcr\noalign{\kern2pt}
$\scriptstyle{#2}$\crcr\noalign{\kern2pt\nointerlineskip}
\kern-2pt$\hfil\displaystyle{#1}\hfil$\crcr}}}\limits}

\begin{document}

\title{$K$-projectors}
\author{Alexander Odesskii}
   \date{}
\vspace{-20mm}
   \maketitle
\vspace{-7mm}
\begin{center}
Brock University, 1812 Sir Isaac Brock Way, St. Catharines, ON, L2S 3A1 Canada\\[1ex]
e-mail: \\
\texttt{aodesski@brocku.ca}
\end{center}

\medskip

\begin{abstract}
We study representations of a free associative algebra $T^*(W\otimes W^*)$ in a vector space $V$ with the property $V\otimes V\cong V\oplus V_0$ where $T^*(W\otimes W^*)$ acts by zero on $V_0$ 
and the tensor product $V\otimes V$ of representations corresponds to the natural homomorphism $W\otimes W^*\to W\otimes W^* \otimes W\otimes W^*$. We develop an algebraic theory of such objects and construct a lot of examples. 

\medskip

\medskip

\end{abstract}

\newpage
\tableofcontents
\newpage

 \section{Introduction}

Let $W$ be a vector space\footnote{All vector spaces in these notes are finite-dimensional and over $\C$.}. Denote by $\A_W$ the category of finite-dimensional representations of $T^*(W\otimes W^*)$, the free associative algebra generated by $W\otimes W^*$, where $W^*$ stands for the dual vector space. An object $V\in Ob(\A_W)$ is just a linear mapping 
$$W\otimes W^*\otimes V\to V$$
without any constrains. By a well-known result from linear algebra it can be considered as a linear map
$$W\otimes V\to V\otimes W$$
as well. The category $\A_W$ possesses a structure of tensor category: given two representations $P_1:~W\otimes V_1\to V_1\otimes W$ and $P_2:~W\otimes V_2\to V_2\otimes W$ we 
construct its tensor product $V_1\otimes V_2$ as the following composition
$$W\otimes V_1\otimes V_2\xrightarrow{P_1\otimes 1} V_1\otimes W\otimes V_2\xrightarrow{1\otimes P_2} V_1\otimes V_2\otimes W.$$
This tensor product often appears in the algebraic study of lattice models \cite{bax}.

We denote by $[V]$ the image of $V\in Ob(\A_W)$ in the Grothendieck group of $\A_W$.

In the paper \cite{K} Maxim Kontsevich posed the following\footnote{See \cite{K}, Section 3.5, Question 6.} 

{\bf Question.} Are there interesting objects $V\in Ob(\A_W)$ for some $W$ such that $[V\otimes V] = [V]$?

In these notes we make an attempt to answer this question.

{\bf Definition 1.} A $K$-projector $P$ is a linear mapping $P:~W\otimes V\to V\otimes W$ such that in the tensor category $\A_W$ defined above we have
$$[V\otimes V] = [V].$$

To create an algebraic theory of $K$-projectors we need to add some constrains to this definition.

{\bf Definition 2.} A $K$-projector $P:~W\otimes V\to V\otimes W$ is called {\em regular} if $V$ is an irreducible representation of $T^*(W\otimes W^*)$ and $V\otimes V \cong V\oplus V_0$ where $[V_0]=0$ in the Grothendieck group of $\A_W$.

It seems that this constrain is not very restrictive: any ``interesting'' $K$-projector should be regular. However, we will also need a more restrictive

{\bf Definition 3.} A $K$-projector $P:~W\otimes V\to V\otimes W$ is called {\em perfect} if $V$ is an irreducible representation of $T^*(W\otimes W^*)$ and $V\otimes V \cong V\oplus V_0$ where $T^*(W\otimes W^*)$ acts by zero on $V_0$.

\section{K-projectors as representations of certain associative algebras}

In order to formulate the main results of this Section we need the following

{\bf Definition 4.} A $P$-algebra is a vector space $V$ together with two linear maps $\mu:~V\otimes V\to V$ and $\lambda:~ V\to V\otimes V$ subject to the following constrains:

- $\mu$ is an associative multiplication and $\lambda$ is an associative comultiplication in $V$.

- The following composition $V\xrightarrow{\lambda}V\otimes V\xrightarrow{\mu}V$ is equal to the identity.

{\bf Remark 1.} In general a $P$-algebra is not a bialgebra. We do not require for $\lambda$ to be a homomorphism of associative algebras.

{\bf Definition 5.} Let $V$ be a $P$-algebra. A weak envelope of $V$ is an associative algebra $En_w(V)$ generated by $V\otimes V^*$ with the following defining relations
\begin{equation}\label{wen1}
\sum_{\alpha} (u\otimes x_{\alpha}^{\prime\prime})( v\otimes x_{\alpha}^{\prime})=\mu(v\otimes u)\otimes x,~\text{where}~\mu^*(x)=\sum_{\alpha} x_{\alpha}^{\prime}\otimes x_{\alpha}^{\prime\prime}
\end{equation}
\begin{equation}\label{wen2}
\sum_{\alpha} (u_{\alpha}^{\prime\prime}\otimes x) (u_{\alpha}^{\prime} \otimes y)=u\otimes \lambda^*(y\otimes x),~\text{where}~\lambda(u)=\sum_{\alpha} u_{\alpha}^{\prime}\otimes u_{\alpha}^{\prime\prime}
\end{equation}
for all $u,v\in V$, $x,y\in V^*$. Here $\mu^*$, $\lambda^*$ stand for the dual maps.

{\bf Definition 6.} Let $V$ be a $P$-algebra. An envelope of $V$ is the tensor product of associative algebras $V$ with multiplication $\bar{\mu}$ and $V^*$ with multiplication $\bar{\lambda^*}$ 
where $\bar{\mu},~\bar{\lambda^*}$ stand for opposite multiplications. In other 
words, an envelope is the vector space $En(V)=V\otimes V^*$ with the associative product defined by  
\begin{equation}\label{en}
(u\otimes x)(v\otimes y)=\mu(v\otimes u)\otimes \lambda ^*(y\otimes x).
 \end{equation}

Let us choose a basis $e_1,...,e_n$ in $V$. Let $\mu(e_i\otimes e_j)=c^k_{ij}e_k$ and $\lambda(e_i)=s_i^{jk}e_j\otimes e_k$. These tensors must satisfy the constrains\footnote{Here and in the sequel we assume summation by repeated latin indexes.}
\begin{equation}\label{rel}
c_{ij}^rc_{rk}^t=c_{ir}^tc_{jk}^r,~~~s_t^{ir}s_r^{jk}=s_r^{ij}s_t^{rk},~~~s_i^{rt}c_{rt}^j=\delta_i^j 
\end{equation}
by Definition 4 where $\delta_i^j$ is a Kronecker delta. Let $e_i^j=e_i\otimes e^j$ be the corresponding basis of $V\otimes V^*$. The weak envelope $En_w(V)$ is generated by $e_i^j,~i,j=1,...,n$ 
with defining relations 
\begin{equation}\label{wenb}
c_{rt}^ke_j^te_i^r=c_{ij}^pe_p^k,~~~s_i^{rt}e_t^ke_r^j=s_p^{jk}e_i^p. 
\end{equation}

The envelope $En(V)$ is the associative algebra with the basis $e_i^j,~i,j=1,...,n$ and the product given by 
\begin{equation}\label{enb}
 e_j^le_i^k=c_{ij}^ts_r^{kl}e_t^r.
\end{equation}

{\bf Lemma 1.} There exists a homomorphism of associative algebras $\phi: ~En_w(V)\to En(V)$ such that $\phi |_{V\otimes V^*}=\text{Id}$. In other words, $En(V)$ is a quotient of $En_w(V)$.

{\bf Proof.} We need to prove that the relations (\ref{wen1}), (\ref{wen2}) are both consequences of (\ref{en}). Computing the products in the l.h.s. of (\ref{wen1}) by virtue of (\ref{en}) we obtain
$$\sum_{\alpha} (u\otimes x_{\alpha}^{\prime})( v\otimes x_{\alpha}^{\prime\prime})=\sum_{\alpha}\mu(v\otimes u)\otimes \lambda^*(x_{\alpha}^{\prime\prime}\otimes x_{\alpha}^{\prime})=\mu(v\otimes u)\otimes \lambda^*(\mu^*(x))=\mu(v\otimes u)\otimes x$$
where we used the property of $P$-algebras $\mu\circ\lambda=\text{Id}$ from the Definition 4. The proof of (\ref{wen2}) is similar. $\Box$

The main results of this section are the following

{\bf Theorem 1.} There exists a natural correspondence: given a {\em regular} $K$-projector $P:~W\otimes V\to V\otimes W$ one can construct a $P$-algebra structure in the vector space $V$ and a representation of $En_w(V)$ in the vector space $W$. Moreover, given a $P$-algebra structure in a vector space $V$ and a representation of $En_w(V)$ in a vector space $W$ obtained from some regular $K$-projector, one can reconstruct this $K$-projector. On the other hand, given an arbitrary $P$-algebra structure in a vector space $V$ and a representation of $En_w(V)$ in a vector space $W$ one can construct a representation of $T^*(W\otimes W^*)$ in the vector space $V$ such that $V\otimes V\cong V\oplus V_0$ but it is not necessarily that $[V_0]=0$ nor $V$ is irreducible as a representation of  $T^*(W\otimes W^*)$.

{\bf Theorem 2.} There exists a natural correspondence: given a {\em perfect} $K$-projector $P:~W\otimes V\to V\otimes W$ one can construct a $P$-algebra structure in the vector space $V$ and a representation of $En(V)$ in the vector space $W$. Moreover, given a $P$-algebra structure in a vector space $V$ and a representation of $En(V)$ in a vector space $W$ obtained from some perfect $K$-projector, one can reconstruct this $K$-projector. On the other hand, given an arbitrary $P$-algebra structure in a vector space $V$ and a representation of $En(V)$ in a vector space $W$ one can construct a representation of $T^*(W\otimes W^*)$ in the vector space $V$ such that $V\otimes V\cong V\oplus V_0$ where $V_0$ is a zero representation of  $T^*(W\otimes W^*)$ but  $V$ is not necessarily irreducible.

{\bf Remark 2.} Therefore, given an arbitrary $P$-algebra structure in the vector space $V$ and a representation of $En(V)$ in a vector space $W$ one can construct a $K$-projector which is not necessarily perfect.

{\bf Proof.} Let $P:~W\otimes V\to V\otimes W$ be a {\em regular} $K$-projector. We have $V\otimes V\cong V\oplus V_0$ and therefore 
$V\otimes V\otimes V\cong V\otimes V\oplus V_0\otimes V\cong V\oplus V_0^{\prime}$ where $[V_0]=[V_0^{\prime}]=0$ and $V$ is an irreducible representation of $T^*(W\otimes W^*)$. Therefore, the both representations  $V\otimes V$ and $V\otimes V\otimes V$ contain $V$ with multiplicity one.  Let $\mu:~V\otimes V\to V$ and $\lambda:~ V\to V\otimes V$ be intertwining operators of these $T^*(W\otimes W^*)$-modules. By the Schur lemma these exist and unique up to multiplication by scalars. Moreover, consider the diagrams
$$V\xrightarrow{\lambda}V\otimes V\xrightarrow{\mu}V$$

$$\begin{array}{ccc}
V\otimes V\otimes V & \xrightarrow{~\mu\otimes 1~} & V\otimes V\\
1\otimes\mu~\Big\downarrow ~~~~~~~~~&  & ~~~\Big\downarrow~\mu  \\
 V\otimes V & \xrightarrow{~~~\mu~~~} &  V   \end{array}$$

$$\begin{array}{ccc}
V & \xrightarrow{~~~\lambda~~~} & V\otimes V\\
\lambda~\Big\downarrow ~~~&  & ~~~~~~~~\Big\downarrow~\lambda\otimes 1  \\
 V\otimes V & \xrightarrow{~~1\otimes\lambda~~} &  V\otimes V\otimes V   \end{array}$$

By the Schur lemma the composition in the first diagram is proportional to identity and can be made equal to identity by rescaling of $\mu$ or $\lambda$. Also by the Schur lemma the second and the third diagrams are commutative up to proportionality but it is easy to see (considering the corresponding pentagon diagrams) that these coefficients of proportionality are equal to one. 
Therefore, we obtain a $P$-algebra. To obtain a representation of the weak envelope $En_w(V)$ we consider our $K$-projector $P:~W\otimes V\to V\otimes W$ as $W\otimes V\otimes V^*\to  W$. This gives a representation of $T^*(V\otimes V^*)$ in $W$. To prove the relations (\ref{wen1}) and (\ref{wen2}) we consider the following diagrams

\begin{equation}\label{d1}
\begin{array}{cccccc}
W\otimes V\otimes V & \xrightarrow{~P\otimes 1~} & V\otimes W\otimes V  &\xrightarrow{~1\otimes P~}&V\otimes V\otimes W\\
1\otimes\mu~\Big\downarrow ~~~~~~~~& & & & ~~~~~~~\Big\downarrow~\mu\otimes 1  \\
 W\otimes V && \xrightarrow{~~~~~~~~~~~P~~~~~~~~~~~} &&  V\otimes W   \end{array} 
\end{equation}

and

\begin{equation}\label{d2}
\begin{array}{cccccc}
W\otimes V\otimes V & \xrightarrow{~P\otimes 1~} & V\otimes W\otimes V  &\xrightarrow{~1\otimes P~}&V\otimes V\otimes W\\
1\otimes\lambda~\Big\uparrow ~~~~~~~~& & & & ~~~~~~~\Big\uparrow~\lambda\otimes 1  \\
 W\otimes V && \xrightarrow{~~~~~~~~~~~P~~~~~~~~~~~} &&  V\otimes W   \end{array} 
\end{equation}

These diagrams are commutative by construction and give the relations (\ref{wen1}) and (\ref{wen2}). In a basis $e_1,...,e_n$ the structure of $T^*(V\otimes V^*)$-module in $W$ is defined by
$$w\otimes e_i\overset{P}\longmapsto e_j\otimes e_i^jw$$
where $e_i^j=e_i\otimes e^j$ is the corresponding basis of $V\otimes V^*$ and $w\in W$. The commutativity of the diagram (\ref{d1}) gives
$$\begin{array}{cccccc}
w\otimes e_i\otimes e_j & \xmapsto{~P\otimes 1~} & e_t\otimes e_i^tw\otimes e_j  &\xmapsto{~1\otimes P~}&e_t\otimes e_r \otimes e_j^re_i^tw\\
\begin{rotate}{-90}$\longmapsto$ \end{rotate}&&&&\begin{rotate}{-90}$\longmapsto$ \end{rotate}\\
1\otimes\mu ~~~~~~~~~~& & & & ~~~~~~~~~~~~~\mu\otimes 1  \\
 c_{ij}^rw\otimes e_r && \xmapsto{~~~~~~~~~~~P~~~~~~} &&  c_{ij}^re_k\otimes e_r^kw=c_{tr}^k e_k \otimes e_j^re_i^tw ~~~~~~~~~~~~~~~  \end{array}  $$
which gives the first relation in (\ref{wenb}). Similar computation with the diagram (\ref{d2}) gives the second relation in (\ref{wenb}).

Therefore, we obtain a representation of $En_w(V)$. It is clear that some of these steps can be done backwards: given a $P$-algebra $V$ and a representation of $En_w(V)$ in a vector space $W$ we obtain the commutative diagrams (\ref{d1}) and (\ref{d2}) which show that $V\otimes V\cong V\oplus V_0$ as $T^*(W\otimes W^*)$-modules. But we do not obtain in general that $[V_0]=0$ nor that $V$ is irreducible. This proves the Theorem 1. $\Box$

Let $P:~W\otimes V\to V\otimes W$ be a {\em perfect} $K$-projector. It is clear that any perfect $K$-projector is also regular so we have a $P$-algebra $V$ and a representation of $En_w(V)$ in $W$ by 
the previous consideration. We need to prove the relations (\ref{en}). It is clear that $\lambda\circ\mu$ is a projector in $V\otimes V$. Indeed, $(\lambda\circ\mu)^2=\lambda\circ(\mu\circ\lambda)\circ\mu=\lambda\circ\text{Id}\circ\mu=\lambda\circ\mu.$ This projector is a homomorphism of $T^*(W\otimes W^*)$-modules and $\text{Im}(\lambda\circ\mu)=V$. Therefore, $\text{Ker}(\lambda\circ\mu)=\text{Im}(\text{Id}-\lambda\circ\mu)=V_0$ which is a zero representation of $T^*(W\otimes W^*)$ by definition of a perfect $K$-projector. 
This means that the composition of the linear maps $P\otimes 1$ and $1\otimes P$ applied to $\text{Im}(\text{Id}-\lambda\circ\mu)$ gives zero which is equivalent to (\ref{en}). 

Computing in a basis $e_1,...,e_n$ of $V$ we get $$(\text{Id}-\lambda\circ\mu)e_i\otimes e_j=e_i\otimes e_j-c_{ij}^rs_r^{pt}e_p\otimes e_t$$ and 
$$w\otimes(e_i\otimes e_j-c_{ij}^rs_r^{pt}e_p\otimes e_t)\overset{P\otimes 1}\longmapsto e_k\otimes e_i^kw\otimes e_j-c_{ij}^rs_r^{pt}e_k\otimes e_p^kw\otimes e_t\overset{1\otimes P}\longmapsto$$
$$\overset{1\otimes P}\longmapsto e_k\otimes e_l\otimes e_j^le_i^kw-c_{ij}^rs_r^{pt}e_k\otimes e_l\otimes e_t^le_p^kw=e_k\otimes e_l\otimes (e_j^le_i^k-c_{ij}^rs_r^{pt}e_t^le_p^k)w=0$$
and we obtain 
\begin{equation}\label{m}
e_j^le_i^k=c_{ij}^rs_r^{pt}e_t^le_p^k 
\end{equation}
Applying the second relation from (\ref{wenb}) to the r.h.s. of (\ref{m})  we obtain (\ref{enb}).
Therefore, we obtain a representation of $En(V)$.

It is clear that some of these steps can be done backward:  given a representation of $En(V)$ in a vector space $W$ one can construct a $K$-projector. The only problem is that $V$ is not irreducible in general and therefore the corresponding $K$-projector is not necessarily perfect. $\Box$

According to \cite{K}, any $K$-projector $P:~W\otimes V\to V\otimes W$ gives rise an infinite series of projectors $P_N:~W^{\otimes N}\to W^{\otimes N}$. In our framework this can be reformulated as follows

{\bf Lemma 2.} Let $V$ be a $P$-algebra. Let $P_1=\text{Id}\in En(V)= V\otimes V^*=\text{End}(V)$ and in general 
$P_N=\text{Id}\in En(V)^{\otimes N}= \text{End}(V^{\otimes N})$. Then  $P_N^2=P_N$ for all $N=1,2,...$. In a basis $e_1,...,e_n$ of $V$ we have 
\begin{equation}\label{pr}
P_1= e^i_i,~P_2=e_i^j\otimes e_j^i~,...,~P_N= e_{i_1}^{i_2}\otimes e_{i_2}^{i_3}\otimes...\otimes e_{i_{N-1}}^{i_N}\otimes e_{i_N}^{i_1},... 
\end{equation}
Therefore, if $W$ is a $En(V)$-module, then $P_N$ defines a projector in $W^{\otimes N}$ for each $N=1,2,...$

{\bf Remark 3.} According to \cite{K}, if $V,~W$ correspond to a $K$-projector $P:~W\otimes V\to V\otimes W$, then the numbers $\text{tr}(P_N)=\text{dim}(\text{Im}(P_N))$ are important invariants of $P$.

{\bf Proof.} We have $$P_1^2=e_{i}^{i}e_{j}^{j}= c_{ij}^rs^{ij}_te^t_r= \delta_t^re_r^t=e_r^r=P_1$$
where we use (\ref{enb}) and (\ref{rel}). The computation for $P_N$ is similar. $\Box$

\section{Examples  of P-algebras}

The results of the previous Section suggest the following way for constructing examples of perfect $K$-projectors:

{\bf 1.} Construct examples of $P$-algebras.

{\bf 2.} Construct representations of the envelope of each $P$-algebra obtained in the previous step.

{\bf 3.} For each representation obtained in the previous step check if the corresponding $K$-projector is perfect.

In this Section we address the first step which seems to be the hardest.

{\bf Example 1.} Let $V$ be a two-dimensional vector space with a basis $e_1,~e_2$. Assume that 
\begin{equation}\label{cop}
 \lambda(e_1)=e_1\otimes e_1,~~~\lambda(e_2)=e_2\otimes e_2.
\end{equation}
Then each of the following formulas defines a $P$-algebra

\begin{equation}\label{p1}
\mu(e_1\otimes e_1)=e_1,~~~\mu(e_2\otimes e_2)=e_2,~~~\mu(e_1\otimes e_2)=\mu(e_2\otimes e_1)=0, 
\end{equation}

\begin{equation}\label{p2}
\mu(e_1\otimes e_1)=e_1,~~~\mu(e_2\otimes e_2)=\mu(e_1\otimes e_2)=\mu(e_2\otimes e_1)=e_2, 
\end{equation}

\begin{equation}\label{p3}
\mu(e_1\otimes e_1)=\mu(e_1\otimes e_2)=e_1,~~~\mu(e_2\otimes e_2)=\mu(e_2\otimes e_1)=e_2, 
\end{equation}

\begin{equation}\label{p4}
\mu(e_1\otimes e_1)=\mu(e_2\otimes e_1)=e_1,~~~\mu(e_2\otimes e_2)=\mu(e_1\otimes e_2)=e_2, 
\end{equation}

{\bf Lemma 3.} All $P$-algebras with $\text{dim}~V=2$ and semisimple comultiplication are listed in Example 1.

{\bf Proof.} Any semisimple two-dimensional algebra is commutative and the coproduct of the corresponding coalgebra is given by (\ref{cop}). One can check by direct computation that any product compatible with the coproduct given by  (\ref{cop}) in the sense of (\ref{rel}) is given by one of (\ref{p1}) - (\ref{p4}). $\Box$

It seems to be natural, at least from the pure algebraic point of view, to try to classify all $P$-algebras such that both product and coproduct are semisimple. The first step in this direction is the following 

{\bf Theorem 3.} Let $V$ be a $P$-algebra such that the product $\mu$ and the coproduct $\lambda$ are both semisimple and commutative. Then there exist bases $e_1,...,e_n$ and $f_1,...,f_n$ of $V$ and a nonsingular $n\times n$ matrix $(r_{\alpha}^{\beta})$ such that 
\begin{equation}\label{ss}
r_{\alpha}^{\beta}\in \{0,~1\}, ~~~e_{\alpha}=r_{\alpha}^tf_t, ~~~\lambda (e_{\alpha})=e_{\alpha}\otimes e_{\alpha},~~~\mu(f_{\beta}\otimes f_{\beta})=f_{\beta}~~~\alpha,\beta=1,...,n.  
\end{equation}
The matrix $(r_{\alpha}^{\beta})$ is determined by a $P$-algebra $V$ up to independent permutations of rows and columns. Moreover, for any nonsingular $(0,1)$-matrix $(r_{\alpha}^{\beta})$ the formulas (\ref{ss}) define a $P$-algebra with commutative semisimple product and coproduct.

{\bf Remark 4.} For $n=2$ there are two nonsingular $(0,1)$-matrices up to independent rows and columns permutations: 
$\left(\begin{array}{cc}                                                                                     
 1 & 0\\0 & 1 
\end{array}\right)$ and 
$\left(\begin{array}{cc}                                                                                     
 1 & 1\\0 & 1 
\end{array}\right)$. The corresponding $P$-algebras are given in Example 1 by (\ref{p1}) and (\ref{p2}).

{\bf Remark 5.} For $n=3$ there are $8$ nonsingular $(0,1)$-matrices up to independent rows and columns permutations.

{\bf Proof.} Let $V$ be a $P$ algebra such that the product $\mu$ and the coproduct $\lambda$ are both semisimple and commutative. It is known \cite{vdw} that there exist bases $e_1,...,e_n$ and $f_1,...,f_n$ of $V$ such that $\lambda (e_{\alpha})=e_{\alpha}\otimes e_{\alpha},~\mu(f_{\beta}\otimes f_{\beta})=f_{\beta},~\alpha,\beta=1,..,n$ and this bases are defined uniquely up to permutations. Define a matrix $(r_{\alpha}^{\beta})$ by $e_{\alpha}=r_{\alpha}^tf_t$, this matrix is nonsingular by construction and is defined uniquely up to independent rows and columns permutations. The only condition we need to require is $\mu\circ\lambda=\text{Id}$. We have 
$$\mu(\lambda(e_{\alpha}))=\mu(e_{\alpha}\otimes e_{\alpha})=r_{\alpha}^pr_{\alpha}^q\mu(f_p\otimes f_q)=\delta_{p,q}r_{\alpha}^pr_{\alpha}^qf_q=(r_{\alpha}^p)^2f_p,~~~e_{\alpha}=r_{\alpha}^pf_p$$
and therefore $(r_{\alpha}^p)^2=r_{\alpha}^p$ which is equivalent to $r_{\alpha}^p\in\{0,1\}$. $\Box$

{\bf Example 2.} Let $V$ be a $P$-algebra such that $V\cong \text{Mat}_2(\C)$ as an algebra with the product $\mu$ and $V\cong \C\oplus\C\oplus\C\oplus\C$ as a coalgebra with the coproduct $\lambda$. Let  $e_1,e_2,e_3,e_4\in\text{Mat}_2(\C)$ be a basis of $\text{Mat}_2(\C)$ such that $\lambda(e_{\alpha})=e_{\alpha}\otimes e_{\alpha},~\alpha=1,2,3,4$. These matrices are defined uniquely up to permutations and simultaneous conjugations. The only condition we have is $e_{\alpha}^2=e_{\alpha},~\alpha=1,2,3,4$. There are two possibilities (up to permutations of $e_{\alpha}$): $\text{rk}(e_1)=2,~\text{rk}(e_2)=\text{rk}(e_3)=\text{rk}(e_4)=1$ and $\text{rk}(e_1)=\text{rk}(e_2)=\text{rk}(e_3)=\text{rk}(e_4)=1$. In the first case we have 
\begin{equation}\label{c1}
e_1=\left(\begin{array}{cc}                                                                                     
 1 & 0\\0 & 1 
\end{array}\right),~e_2=\left(\begin{array}{cc}                                                                                     
 1 & 0\\0 & 0 
\end{array}\right),~e_3=\left(\begin{array}{cc}                                                                                     
 a_1 & b_1\\c_1 & 1-a_1 
\end{array}\right),~e_4=\left(\begin{array}{cc}                                                                                     
 a_2 & b_2\\c_2 & 1-a_2 
\end{array}\right)
\end{equation}
where $a_1(1-a_1)=b_1c_1,~a_2(1-a_2)=b_2c_2$ and in the second case we have 
\begin{equation}\label{c2}
e_1=\left(\begin{array}{cc}                                                                                     
 1 & 0\\0 & 0 
\end{array}\right),~e_2=\left(\begin{array}{cc}                                                                                     
 a_1 & b_1\\c_1 & 1-a_1 
\end{array}\right),~e_3=\left(\begin{array}{cc}                                                                                     
 a_2 & b_2\\c_2 & 1-a_2 
\end{array}\right),~e_4=\left(\begin{array}{cc}                                                                                     
 a_3 & b_3\\c_3 & 1-a_3 
\end{array}\right)
\end{equation}
where $a_1(1-a_1)=b_1c_1,~a_2(1-a_2)=b_2c_2,~a_3(1-a_3)=b_3c_3$.

{\bf Example 3.} Similarly to the previous example, let $V$ be a $P$-algebra such that $V\cong \text{Mat}_n(\C)$ as an algebra with the product $\mu$ and let $V$ be a commutative semisimple coalgebra with respect to the coproduct $\lambda$. It is clear that the description of such $P$-algebras is the same as the description of bases in $\text{Mat}_n(\C)$ consisting of projectors up to permutation and conjugation.

{\bf Example 4.} Let $V\cong \text{Mat}_n(\C)$ as an algebra with the product $\mu$ and as a coalgebra with the coproduct $\lambda$. To describe such $P$-algebras we need to describe bases $e_{ij},~i,j=1,...,n$ in $\text{Mat}_n(\C)$ such that $e_{i1}e_{1j}+...+e_{in}e_{nj}=e_{ij},~i,j=1,...,n$. In other words, we need to describe block $n\times n$-matrices $E=(e_{ij})$ consisting of $n\times n$-matrices $e_{ij}$ such that $E^2=E$ and $e_{ij}$ is a basis of $\text{Mat}_n(\C)$. Two block matrices $E$ and $E^{\prime}$ define isomorphic $P$-algebras iff $E^{\prime}=UEU^{-1}$ where $U$ is a product of a block diagonal matrix with equal matrices on the diagonal and a block matrix with scalar entries.

{\bf Example 5.} Set $n=2$ in the previous example. We have $E=\left(\begin{array}{cc}                                                                                     
 a & b\\c & d 
\end{array}\right)$ where $a,b,c,d\in \text{Mat}_2(\C)$ is a basis. Assume that $b$ is invertible. In this case the condition $E^2=E$ is equivalent to $d=1-b^{-1}ab,~c=b^{-1}a(1-a)$ for arbitrary $a,~b$. In this case $\text{tr}(E)=\text{tr}(a)+\text{tr}(d)=2$ and therefore $\text{rk}(E)=2$.

{\bf Remark 6.} It is clear from these examples that the problem of classification of $P$-algebras with semisimple product and coproduct is wild. However, it is easy to construct examples of such $P$-algebras.

{\bf Theorem 4.} There exists one-to-one correspondence between $P$-algebras with simisimple products and coproducts and the following data:

- Two families of vector spaces $L_1,...,L_r$ and $M_1,...,M_t$ such that 
$$\sum_{\alpha=1}^r(\text{dim} L_{\alpha})^2=\sum_{\beta=1}^t(\text{dim} M_{\beta})^2$$

- An isomorphism of vector spaces
$$Q:~\bigoplus_{\alpha=1}^rL_{\alpha}\otimes L_{\alpha}^*\to\bigoplus_{\beta=1}^tM_{\beta}\otimes M_{\beta}^*$$
such that all induced homomorphisms  $Q_{\alpha,\beta}:~L_{\alpha}\otimes L_{\alpha}^*\to M_{\beta}\otimes M_{\beta}^*$ considered as linear maps $\bar{Q}_{\alpha,\beta}:~L_{\alpha}\otimes M_{\beta}^*\to L_{\alpha}\otimes M_{\beta}^*$ are projectors\footnote{We agree that zero maps are also projectors.}: $\bar{Q}_{\alpha,\beta}^2=\bar{Q}_{\alpha,\beta}$.

{\bf Proof.} We are given that $V$ is a semisimple associative algebra with respect to the product $\mu$. By the Wedderburn's theorem \cite{vdw} we have an isomorphism of associative algebras 
$$V\cong \bigoplus_{\beta=1}^t\text{Hom}(M_{\beta},M_{\beta})$$
where $M_1,...,M_t$ are some vector spaces. Similarly we have
$$V^*\cong \bigoplus_{\alpha=1}^r\text{Hom}(L_{\alpha},L_{\alpha})$$
because the algebra $V^*$ with the product $\lambda^*$ is also semisimple. Using canonical isomorphisms $\text{Hom}(U,U)\cong U\otimes U^*$ and $(U\otimes U^*)^*\cong U\otimes U^*$ we can rewrite these isomorphisms as
$$V\cong \bigoplus_{\alpha=1}^rL_{\alpha}\otimes L_{\alpha}^*~~~\text{and}~~~V\cong \bigoplus_{\beta=1}^tM_{\beta}\otimes M_{\beta}^*$$
which gives the isomorphism $Q$. The only property of these data we need to require is $\mu\circ\lambda=\text{Id}$. Let $e_{1,\alpha},...,e_{k_{\alpha},\alpha}$ be a basis of $L_{\alpha}$ and let 
$f_{1,\beta},...,f_{m_{\beta},\beta}$ be a basis of $M_{\beta}$. We denote by $e_{\alpha}^1,...,e_{\alpha}^{k_{\alpha}}$ and $f_{\beta}^1,...,f_{\beta}^{m_{\beta}}$ the corresponding dual bases. 
We have
$$\mu(\lambda(e_{i,\alpha}\otimes e_{\alpha}^j))=\mu(e_{i,\alpha}\otimes e_{\alpha}^p\otimes e_{p,\alpha}\otimes e_{\alpha}^j)=\sum_{\beta=1}^tQ_{\alpha,i,r}^{\beta,p,t}~Q_{\alpha,p,q}^{\beta,j,s}~m(f_{\beta,t}\otimes f_{\beta}^r\otimes f_{\beta,s}\otimes f_{\beta}^q)=$$$$=\sum_{\beta=1}^tQ_{\alpha,i,r}^{\beta,p,t}~Q_{\alpha,p,q}^{\beta,j,s}~\delta_s^r~f_{\beta,t}\otimes f_{\beta}^q=\sum_{\beta=1}^tQ_{\alpha,i,r}^{\beta,p,t}~Q_{\alpha,p,q}^{\beta,j,r}~f_{\beta,t}\otimes f_{\beta}^q$$
and
$$e_{i,\alpha}\otimes e_{\alpha}^j=\sum_{\beta=1}^tQ_{\alpha,i,q}^{\beta,j,t}~f_{\beta,t}\otimes f_{\beta}^q.$$
Therefore, we obtain $Q_{\alpha,i,r}^{\beta,p,t}~Q_{\alpha,p,q}^{\beta,j,r}=Q_{\alpha,i,q}^{\beta,j,t}$ where $Q_{\alpha,i,q}^{\beta,j,t}$ are matrix elements of both $Q_{\alpha,\beta}$ and $\bar{Q}_{\alpha,\beta}$. This is equivalent to the statement that $\bar{Q}_{\alpha,\beta}$ are projectors. $\Box$

\section{Examples of K-projectors}

In this Section we construct examples of $K$-projectors using our previous results. We also compute the rank of the projectors $P_N$ in each of these examples\footnote{See Lemma 2.}.

\subsection{The case of commutative semisimple product and coproduct} 

Let $V$ be a $P$-algebra described in Theorem 3. The envelop $En(V)$ is a commutative semisimple algebra with a basis $g_i^j=f_i\otimes e^j$ and a product $g_i^jg_k^l=\delta_{ik}\delta_{jl}g_i^j$. Any representation $W$ of $En(V)$ is a direct sum of one-dimensional representations. Let $W_{ij}$ be a one-dimensional representation of $En(V)$ such that $g_i^j\mapsto 1$ and all other basis elements are mapped to zero. We have 
$$W=\bigoplus_{i,j=1}^n W_{ij}^{\oplus m_{ij}}$$
where $m_{ij}\in\{0,1,2,...\}$ are multiplicities of $W_{ij}$ in $W$. The corresponding $P$-projector is defined by
$$P(w\otimes f_i)=e_j\otimes g_i^jw=r_j^kf_k\otimes g_i^jw.$$
Let $w_{ij}^l,~i,j=1,...,n,~l=1,...,m_{ij}$ be a basis of $W$ such that $w_{ij}^l\in W_{ij}^{\oplus m_{ij}}$. We have 
$$P(w_{jk}^l\otimes  f_i)=\delta_{ij}\delta_{rk}e_r\otimes w_{jk}^l=\delta_{ij}e_k\otimes w_{jk}^l=\delta_{ij}r_k^tf_t\otimes w_{jk}^l.$$
It is clear that the $K$-projector $P$ is perfect iff $V$ is irreducible with respect to the set of operators $S_{jk}:~f_i\mapsto \delta_{ij}r_k^tf_t,~j,k=1,...,n$.

Let us compute the projectors $P_N,~N=1,2,...$ defined by (\ref{pr}) in this case. We have 
$$P_1w_{jk}^l=f_i^iw_{jk}^l=\delta_{ij}r_k^iw_{jk}^l=r_k^jw_{jk}^l~~~\text{and}~~~\text{tr}(P_1)=\sum_{j,k=1}^n m_{jk}r_k^j.$$
In general
$$P_Nw_{j_1k_1}^{l_1}\otimes...\otimes w_{j_Nk_N}^{l_N}=f_{i_1}^{i_2}w_{j_1k_1}^{l_1}\otimes...\otimes f_{i_N}^{i_1}w_{j_Nk_N}^{l_N}=\delta_{i_1j_1}r_{k_1}^{i_2}...\delta_{i_Nj_N}r_{k_N}^{i_1}w_{j_1k_1}^{l_1}\otimes...\otimes w_{j_Nk_N}^{l_N}$$
and we obtain
$$P_Nw_{j_1k_1}^{l_1}\otimes...\otimes w_{j_Nk_N}^{l_N}=r_{k_1}^{j_2}r_{k_2}^{j_3}...r_{k_N}^{j_1}w_{j_1k_1}^{l_1}\otimes...\otimes w_{j_Nk_N}^{l_N}.$$
Therefore
\begin{equation}\label{tr1}
 \text{tr}(P_N)=\sum_{j_1,...,k_N=1}^n m_{j_1k_1}r_{k_1}^{j_2}...m_{j_Nk_N}r_{k_N}^{j_1}.
\end{equation}
Let $T$ be an $n\times n$-matrix such that 
$$T_{ij}=\sum_{t=1}^nm_{it}r_t^j.$$
The formula (\ref{tr1}) can be written as
\begin{equation}\label{tr2}
 \text{tr}(P_N)=\text{tr}(T^N).
\end{equation}

\subsection{The general case of semisimple product and coproduct} 

Let $V$ be a $P$-algebra described in Theorem 4. The envelop $En(V)$ is a semisimple algebra 
$$En(V)=V\otimes V^*=\bigoplus_{\substack{1\leq \alpha\leq r\\1\leq\beta\leq t}} L_{\alpha}\otimes L_{\alpha}^*\otimes M_{\beta}\otimes M_{\beta}^*=\bigoplus_{\substack{1\leq \alpha\leq r\\1\leq\beta\leq t}}\text{Hom}(L_{\alpha}\otimes M_{\beta},L_{\alpha}\otimes M_{\beta})$$
with natural multiplication. Any $En(V)$-module is a direct sum of simple modules of the form $L_{\alpha}\otimes M_{\beta}$. Let 
$$W\cong \bigoplus_{\substack{1\leq \alpha\leq r\\1\leq\beta\leq t}} (L_{\alpha}\otimes M_{\beta})^{\oplus m_{\alpha\beta}}$$
where $m_{\alpha\beta}\in\{0,1,2,...\}$ are multiplicities. The corresponding $K$-projector $P:~W\otimes V\to V\otimes W$ written in the form $W\otimes V\otimes V^*\to W$ is the direct sum of 
natural homomorphisms 
$$L_{\alpha}\otimes M_{\beta}\otimes \text{Hom}(L_{\alpha}\otimes M_{\beta},L_{\alpha}\otimes M_{\beta})\to L_{\alpha}\otimes M_{\beta}.$$
It can also be written as the direct sum of compositions
$$(L_{\alpha}\otimes M_{\beta})\otimes (M_{\beta}\otimes M_{\beta}^*)\to (L_{\alpha}\otimes L_{\alpha}^*)\otimes (L_{\alpha}\otimes M_{\beta})\xrightarrow{(\oplus_{\gamma=1}^t Q_{\alpha\gamma})\otimes 1} \left(\bigoplus_{\gamma=1}^t M_{\gamma}\otimes M_{\gamma}^*\right)\otimes (L_{\alpha}\otimes M_{\beta})$$
where the first homomorphism is induced by the composition of  the natural maps $M_{\beta}\otimes M_{\beta}^*\to\C\to L_{\alpha}\otimes L_{\alpha}^*$.

Recall that we denote by $e_{1,\alpha},...,e_{k_{\alpha},\alpha}$ a basis of $L_{\alpha}$ and by 
$f_{1,\beta},...,f_{m_{\beta},\beta}$ a basis of $M_{\beta}$. We denote by $e_{\alpha}^1,...,e_{\alpha}^{k_{\alpha}}$ and $f_{\beta}^1,...,f_{\beta}^{m_{\beta}}$ the corresponding dual bases. In these bases we have
$$P:~(e_{i\alpha}\otimes f_{j\beta})\otimes (e_{k\alpha}\otimes e_{\alpha}^l)\mapsto\delta_i^l~(f_{j\beta}\otimes f_{\beta}^p)\otimes(e_{k\alpha}\otimes f_{p\beta})=\delta_i^l~Q_{\gamma tj}^{\beta qp}(e_{q\gamma}\otimes e_{\gamma}^t)\otimes(e_{k\alpha}\otimes f_{p\beta})$$

The projector $P_1:~W\to W$ defined by (\ref{pr}) is the direct sum of the projectors $\bar{Q}_{\alpha\beta}$, we have 
$$P_1=\bigoplus_{\substack{1\leq \alpha\leq r\\1\leq\beta\leq t}} \bar{Q}_{\alpha\beta}^{\oplus m_{\alpha\beta}}.$$
In particular 
$$\text{tr}(P_1)=\sum_{\substack{1\leq\alpha\leq r\\1\leq\beta\leq t}}m_{\alpha\beta}\text{tr}(\bar{Q}_{\alpha\beta}).$$
More generally, the projector $P_N:~W^{\otimes N}\to W^{\otimes N}$ is given by 
$$P_N=\bigoplus_{\substack{1\leq \alpha_1,...,\alpha_N\leq r\\1\leq\beta_1,...,\beta_N\leq t}} \bar{Q}_{\alpha_1\beta_2}^{\oplus m_{\alpha_1\beta_1}}\otimes...\otimes \bar{Q}_{\alpha_N\beta_1}^{\oplus m_{\alpha_N\beta_N}}$$
and therefore
\begin{equation}\label{tr3}
 \text{tr}(P_N)=\sum_{\substack{1\leq \alpha_1,...,\alpha_N\leq r\\1\leq\beta_1,...,\beta_N\leq t}} m_{\alpha_1\beta_1}\text{tr}(\bar{Q}_{\alpha_1\beta_2})... m_{\alpha_N\beta_N}\text{tr}(\bar{Q}_{\alpha_N\beta_1}).
\end{equation}
Let $T$ be an $r\times r$-matrix such that 
$$T_{\alpha_1\alpha_2}=\sum_{\beta=1}^tm_{\alpha_1\beta}\text{tr}(\bar{Q}_{\alpha_2\beta}).$$
The formula (\ref{tr3}) can be written as
\begin{equation}\label{tr4}
 \text{tr}(P_N)=\text{tr}(T^N).
\end{equation}

\section{Conclusion}

It looks that we have studied the most natural class of $K$-projectors. Indeed, we assumed that $V\otimes V\cong V\oplus {\bf 0}$ and $V$ is a simple $T^*(W\otimes W^*)$-module. Moreover, we assumed that the products $\mu:~ V\otimes V\to V$ and $\lambda^*:~V^*\otimes V^*\to V^*$ are both semisimple. These mean that we have studied ``the most semisimple'' case of $K$-projectors. We have developed a theory which has produced a huge amount of examples of such $K$-projectors, and it is clear that there are so many examples that the classification problem is ``wild''. Indeed, even in the case of semisimple commutative products $\mu$ and $\lambda^*$, wherein the set of $K$-projectors is discrete and is parameterized by nonsingular $(0,1)$-matrices up to independent permutations of rows and columns (along with some additional data, see Theorem 3 and Section 4.1), the combinatorial problem of explicit description (or even counting) of such matrices looks very hard. It is easy however to construct examples of such matrices of arbitrary size, and it is also easy to construct the corresponding examples of $K$-projectors. In the general case of semisimple products $\mu$ and $\lambda^*$, the set of $K$-projectors is parametrized by both discrete and continuous parameters and, as it is demonstrated by Examples 2 to 5, the classification problem is hopeless in general although it is easy to construct examples. We can conclude that the answer to Question 6 from \cite{K} is definitely ``Yes, there are!''.

\section*{Acknowledgments}

I am grateful to Maxim Kontsevich for useful discussions. This work was started when I visited IHES, Bures-sur-Yvette, France. I am grateful to IHES for hospitality and excellent working conditions.


\addcontentsline{toc}{section}{References}

\begin{thebibliography}{9}

\bibitem{bax} R. Baxter, Exactly Solved Models in Statistical Mechanics, Acad. Press, 1981
\bibitem{K} Maxim Kontsevich,  Notes on motives in finite characteristic. Algebra, arithmetic, and geometry: in honor of Yu. I. Manin. Vol. II, 213–247, Progr. Math., 270, Birkhäuser Boston, Inc., Boston, MA, 2009
\bibitem{vdw} van der Waerden, Bartel Leendert, Modern Algebra, New York, N. Y.: Frederick Ungar Publishing Co., 1949


\end{thebibliography}
\end{document}